\documentclass[journal]{IEEEtran}
\usepackage{graphicx}
\usepackage{amsmath}
\usepackage{tabularx}
\usepackage{multirow}
\usepackage{multicol}
\usepackage{lipsum}
\usepackage{color,soul}
\usepackage{flushend}
\allowdisplaybreaks
\begin{document}
\title{Robust Optimal Eco-driving Control\\with Uncertain Traffic Signal Timing}
\author{Chao~Sun,
	Xinwei~Shen,
	and~Scott~Moura,~\IEEEmembership{Member,~IEEE}
\thanks{C. Sun and S. Moura are with the Department
		of Civil and Environmental Engineering, University of California, Berkeley,
		CA, 94704 USA. e-mail: chaosuncal@berkeley.edu, smoura@berkeley.edu.}
\thanks{X. Shen is with the Tsinghua-Berkeley Shenzhen Institute, Shenzhen, China. email: sxw.tbsi@sz.tsinghua.edu.cn.}
\thanks{Manuscript received September 17, 2017.}}

\markboth{\LaTeX\ Class Files, ACC~2018}%
{Shell \MakeLowercase{\textit{et al.}}: Bare Demo of IEEEtran.cls for Journals}

\maketitle

\begin{abstract}
This paper proposes a robust optimal eco-driving control strategy considering multiple signalized intersections with uncertain traffic signal timing. A spatial vehicle velocity profile optimization formulation is developed to minimize the global fuel consumption, with driving time as one state variable. We introduce the concept of ‘effective red-light duration’ (ERD), formulated as a random variable, to describe the feasible passing time through signalized intersections. A chance constraint is appended to the optimal control problem to incorporate robustness with respect to uncertain signal timing. The optimal eco-driving control problem is solved via dynamic programming (DP). Simulation results demonstrate that the optimal eco-driving can save fuel consumption by 50-57\% while maintaining arrival time at the same level, compared with a modified intelligent driver model as the benchmark. The robust formulation significantly reduces traffic intersection violations, in the face of uncertain signal timing, with small sacrifice on fuel economy compared to a non-robust approach.
\end{abstract}
\IEEEpeerreviewmaketitle

\section{Introduction}
Connected and automated vehicle (CAV) technology is revolutionizing the automotive industry. In particular, CAVs may significantly improve safety, energy economy, and convenience. CAVs are able to realize autonomous driving, vehicle to infrastructure (V2I) communication and/or intelligent path/velocity planning \cite{sciarretta2015optimal}. Optimal eco-driving control -- a novel technology brought by CAVs -- is defined as a velocity control method to achieve the most economical fuel, energy or cost performances \cite{jin2016power}. Intuitively speaking, optimal eco-driving seeks the best velocity profile, in some sense, over a specific driving mission. Fig. 1 illustrates the optimal eco-driving concept through V2I communication with a number of traffic signals incorporated. 
\begin{figure}[t]
	\begin{center}
		\includegraphics[trim = 07mm 05mm 07mm 05mm, clip, width=0.49\textwidth]{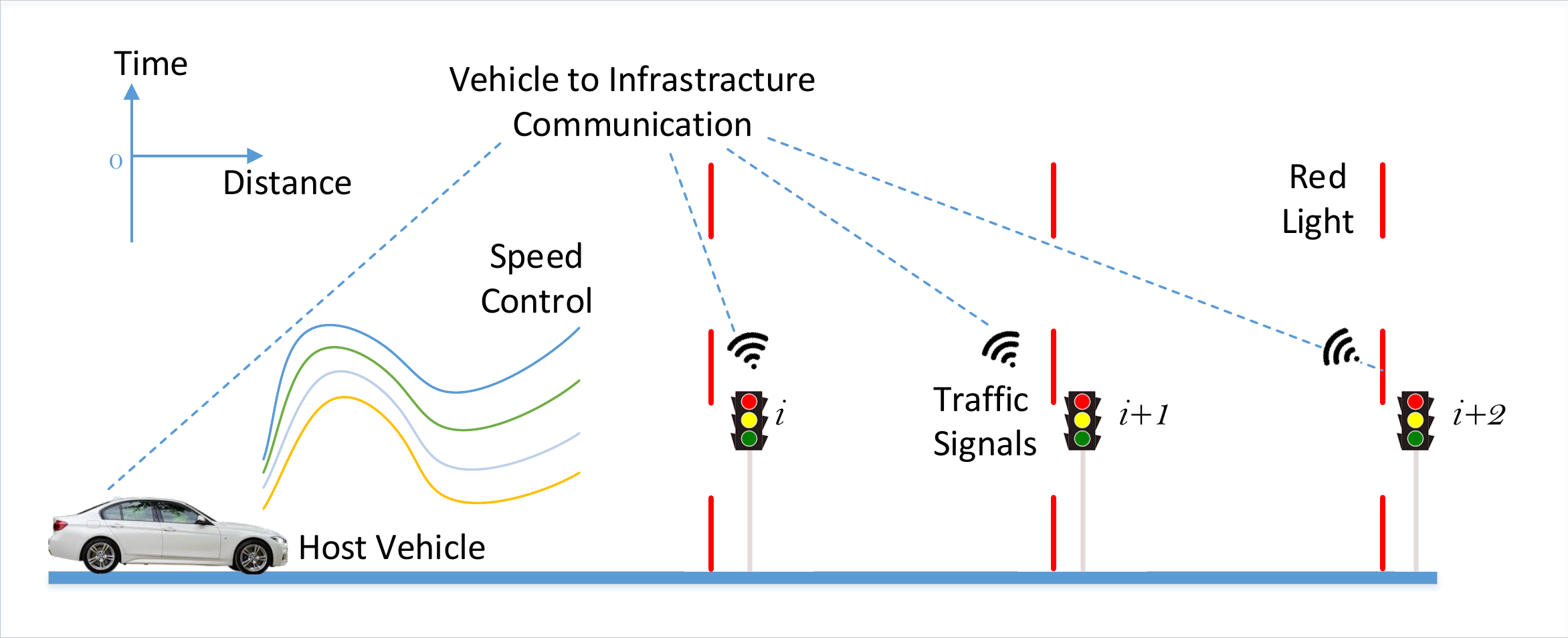}
		\caption{Car optimal eco-driving based on V2I communication with multiple traffic signals incorporated, with vehicle speed control as the main task.}
		\label{fig:ecodriving}
	\end{center}
\end{figure}

In the literature, optimal eco-driving is also known as ecological driving, speed trajectory planning, driving advisory or driver assistance systems. Over the past 10 years, in particular, the optimal eco-driving problem has been intensively studied in the published literature. A heuristic optimal eco-driving strategy is proposed in \cite{saboohi2009model} to minimize the vehicle fuel consumption based on instantaneous fuel performance. With velocity constraints derived from real driving data, another dynamic programming (DP) based optimal eco-driving control is developed in \cite{mensing2013trajectory} for trajectory optimization of an internal combustion engine (ICE) vehicle. Similar approaches are found in \cite{dib2014optimal,kuriyama2010theoretical}, for the optimal energy management as well as speed control of electric vehicles (EV). Experimental results showed a significant increase in energy efficiency. More comprehensively, a cloud-based velocity profile optimization approach is designed in \cite{ozatay2014cloud}, under a spatial domain formulation. Historical velocity data is gathered for speed advising. Spatial domain optimization is further adopted by \cite{lim2017distance} for ecological driving. Uniquely, a short-term adaptation level is added to avoid traffic congestion. Optimal eco-driving has also been integrated into the energy management strategy of hybrid electric vehicles, with interactive Pontryagin’s Minimal Principle used to solve the optimization problem in \cite{hu2016integrated}.

Signal phase and timing (SPaT) information is critical in addressing the optimal eco-driving problem. In \cite{homchaudhuri2015fuel} and \cite{kamal2010board}, hierarchical model predictive control (MPC) is employed for eco-driving in varying traffic environments. Assume the SPaT information is known \textit{a priori}, \cite{li2015eco} solved the optimal eco-departing problem at signalized intersections. Furthermore, \cite{butakov2016personalized} developed a sophisticated on-board driver assistance, which is able to calculate the optimal speed profile with deterministic traffic signals. By considering each signalized intersection as one stage, a multi-stage pseudospectral control method is proposed by \cite{he2015optimal} in an arterial road structure. Hierarchical MPC is also adopted in \cite{homchaudhuri2017fast}, and has demonstrated effective online eco-driving control capabilities. Reference \cite{yang2017eco} considers the car waiting queue in a multi-lane road scenario, and designed an eco-cooperative adaptive cruise control scheme. With a simplified powertrain model and assuming the engine mainly operates along the optimal brake specific fuel consumption (BSFC) line, sequential convex optimization therefore is applied to speed trajectory planning \cite{huang2017speed}.

In the aforementioned studies, signalized intersections are either not considered, or the SPaT information is assumed to be deterministic in the optimal eco-driving control. Ideally, when CAVs have realized V2I communication, SPaT can be communicated to vehicles for optimal eco-driving. This future, however, would require significant penetration of V2I-equipped intersections, which may take decades to realize. Even with V2I-equipped intersections, uncertainty exists due to car waiting queue, pedestrians, bicyclists, varying patterns of traffic lights and other factors, as demonstrated in Fig. 2.  Moreover, an optimal eco-driving approach assuming deterministic SPaT will often pass through intersections exactly at the phase transitions, and thus risks collision. The issue of SPaT uncertainty in optimal eco-routing is significant, and not fully addressed in the existing literature. 
\begin{figure}[t]
	\begin{center}
		\includegraphics[trim = 07mm 07mm 07mm 05mm, clip, width=0.45\textwidth]{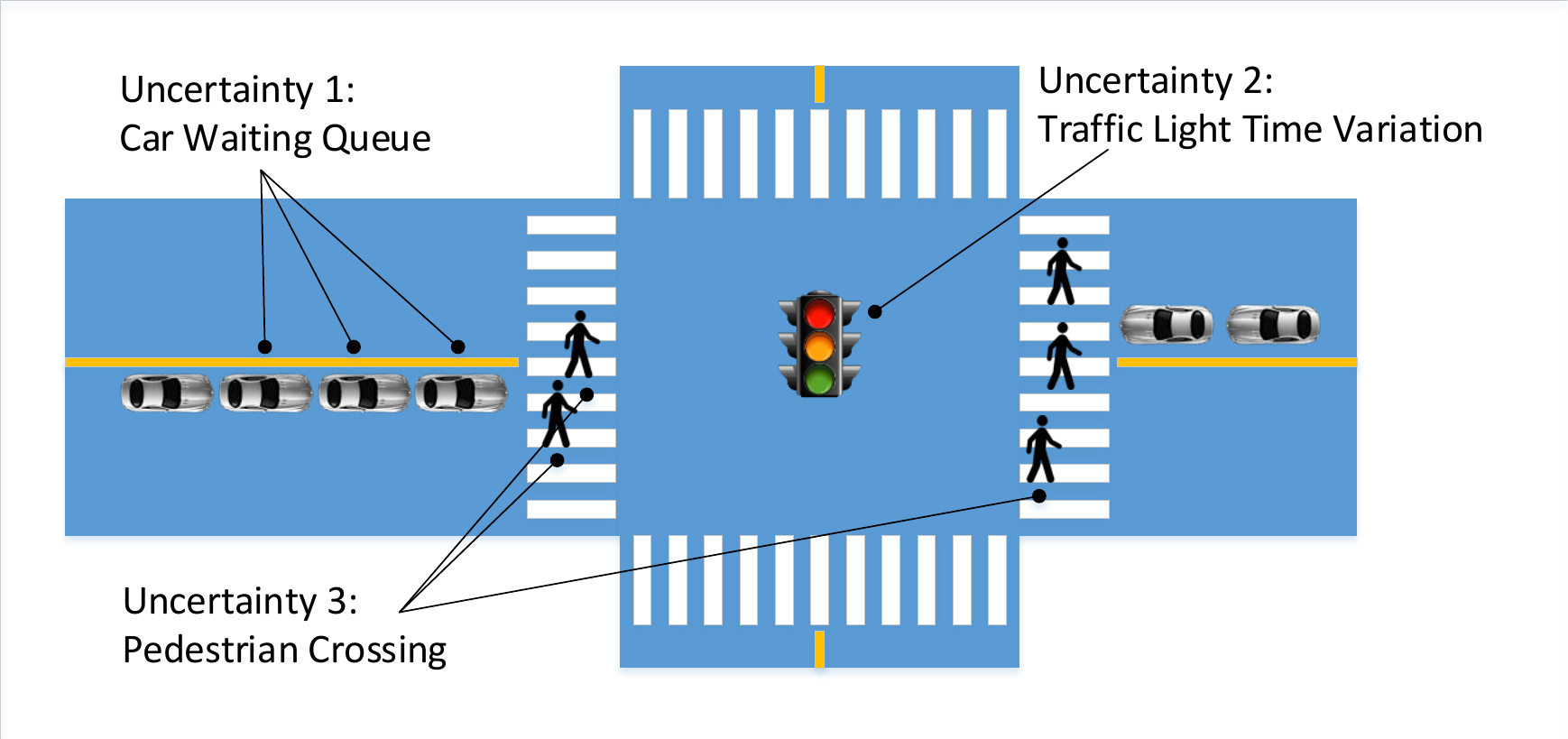}
		\caption{Uncertain factors when passing the road intersections.}
		\label{fig:ERDuncertainty}
	\end{center}
\end{figure}

This paper investigates a fuel-minimizing eco-driving approach that is robust to uncertain feasible vehicle passing times through multiple signalized intersections. The goal is to simultaneously achieve energy economy and safety. The main contributions include:
\begin{itemize} 
	\item
    ‘Effective red-light duration’ (ERD) is proposed to describe the stochastic feasible passing time of vehicles at signalized intersections, composed of a deterministic base red-light duration and a random delay;
	\item
    Signalized intersections are modeled and integrated into the spatial optimal eco-driving formulation, which eliminates the requirement for prior knowledge of accurate arrival time;
	\item
    A robust optimal eco-driving control variant is developed and solved via the dynamic programming. The controller robustness -- and therefore safety -- is significantly improved with little sacrifice of fuel economy.
\end{itemize}
The remainder of the paper is organized as follows. Section II describes the vehicle, traffic signal and driver models. Section III introduces the spatial optimal eco-driving control strategy. Section IV details a robust formulation that considers uncertain feasible passing time at signalized intersections. Section V exhibits the main results, and Section VI draws the main conclusions and future work.

\section{Vehicle, Traffic Signal and Driver Modeling}
\subsection{Vehicle dynamics}
The subject vehicle is equipped with a gasoline ICE and a 6-speed gearbox. Since speed control is the main objective of optimal eco-driving, we consider longitudinal vehicle dynamics and disregard the lateral dynamics. The longitudinal acceleration is calculated by
\begin{equation}
ma = \frac{r_{gb}T_{eng}}{R_{whl}} - mgcos\left(\theta\right)C_{r} - mgsin\left(\theta\right) - \frac{1}{2}\rho A C_{d} v^{2} - T_{brk}
\end{equation}
\begin{equation}
C_{r} = C_{r1} + C_{r2}v
\end{equation}
where $m$ is the vehicle mass, $a$ is the acceleration, $r_{gb}$ is the integrated ratio of gearbox and final drive, $T_{eng}$ is the ICE output torque, $R_{whl}$ is the rolling radius of wheel, $g$ is the gravitational acceleration, $\theta$ is the road grade, and $C_{r}$ is the rolling resistance coefficient. Parameters $\rho$, $A$, $C_{d}$ are the air density, frontal area, and air-dragging resistance coefficient respectively. Variable $v$ is the vehicle velocity, $T_{brk}$ is the braking force enforced on the wheels, $C_{r1}$ and $C_{r2}$ are rolling resistance constants. The longitudinal velocity is computed by
\begin{equation}
v = \frac{\omega_{eng}}{r_{gb}}
\end{equation}
where $\omega_{eng}$ is the ICE rotation speed. The ICE fuel consumption is modeled as a nonlinear map $\phi(\cdot, \cdot)$ that depends on the engine torque and speed:
\begin{equation}
\dot{m}_{fuel} = \psi(T_{eng},\omega_{eng})
\end{equation}
where $\dot{m}_{fuel}$ is the instantaneous fuel consumption, and $\psi$ is the pre-stored fuel map (e.g. a look-up table). The transmission efficiency is ignored in this study. Assume $r_{fd}$ is the final drive ratio. The integrated transmission ratio is formulated as a function of the gear number $N_{gb}$,
\begin{equation}
r_{gb} = f(N_{gb})r_{fd},~~N_{gb}\in \{1,2,3,4,5,6\}
\end{equation}

\subsection{Traffic signal model}
The traffic signal at an intersection is a spatial-temporal system in the optimal eco-driving control problem. Assume the total length of the target driving route is $D_{f}$. The position of the $i$th traffic signal is noted as $D^{i}$ if we treat the signalized intersection as a single point on the road. Therefore,
\begin{equation}
D^{i}\in[0,D_{f}],~~i=\{1,2,3,4,5...I\}
\end{equation}
where $I$ is the total number of traffic signals along the route.

Each traffic signal is modeled with an independent signal-cycling clock in this paper. The universal traveling time of the vehicle is denoted as $t\in R$, and the periodic cycling clock time of the $i$th traffic signal has a period of $c^{i}_{f}\in R$ (clock time zero denotes the beginning of the red light phase). Normally, the period $c^{i}_{f}$ varies at different intersections. The red-light duration is denote by $c^{i}_{r}$. Then we have
\begin{equation}
c^{i}_{r} \in [0,c^{i}_{f}]
\end{equation}

Consider the time when the vehicle departs from its origin. Denote by $c^{i}_{0}$ the periodic signal clock time at this moment. Suppose $t^{i}_{p}$ is the time at which the subject vehicle passes through the $i$th intersection in the universal time domain. We can compute the corresponding time in the periodic traffic signal clock timing by
\begin{equation}
c^{i}_{p} = (c^{i}_{0} + t^{i}_{p})~\text{mod}~c^{i}_{f}
\end{equation}
where $c^{i}_{p}$  is the vehicle passing time in the signal-cycling clock. The modulo operator allows for conversion from the universal time domain to the periodic traffic signal clock time domain. Note that un-signalized intersections or crossings can also be integrated into the model above, which might require on-board cameras or radars to detect the passing conditions.

\subsection{Modified intelligent driver model}
A modified intelligent driver model (IDM) is introduced for comparison with the optimal eco-driving, by imitating human driving behaviors. IDM is originally developed by Treiber \textit{et. al.}, based on the computation of desired distance between the subject vehicle and the vehicle in front or speed limit \cite{kesting2010enhanced}. We enhanced the driver model with an ability to preview traffic signals and adjust speed accordingly. Assume the desired distance between the subject vehicle and front vehicle is $D_{des}$, then
\begin{equation}
D_{des} = D^{min}_{des} + v \cdot t_{hw}-\frac{v D_{sf}}{2\sqrt{a^{max} a_{c}}}
\end{equation}
where $D^{min}_{des}$ is the minimal vehicle distance, $t_{hw}$ is the desired time headway to the preceding vehicle, $D_{sf}$ is the real distance between the subject vehicle and preceding vehicle, $a^{max}$ is the maximal vehicle acceleration ability, and $a_{c}$ is the preferred deceleration for comfort.

The vehicle acceleration at each time step is computed by comparing the desired gap distance with the current distance. An additional speed limit term is added to ensure safety,
\begin{equation}  \label{eq:old_idm}
a = a^{max} \left[1 - \left(\frac{v}{v^{max}}\right)^4 - \left(\frac{D_{des}}{D_{sf}}\right)^2 \right]
\end{equation}

To interact with traffic signals or stop signs, we modify the IDM by enabling the driver model to preview the traffic signal or stop line status at a human-vision distance $D_{v}$. Assume the current location of the vehicle is $D$, the vehicle longitudinal velocity dynamics in (\ref{eq:old_idm}) thus becomes
\begin{equation}
 a = \begin{cases}
a^{max} [1 - \left(\frac{v}{v^{max}}\right)^4 - \left(\frac{D_{des}}{D_{sf}}\right)^2], & \text{if $S_{tss}(D+D_{v})=0$}.\\
- \frac{v^2}{2 D_{sf}}, & \text{if $S_{tss}(D+D_{v})=1$}.
\end{cases}
\end{equation}
where $S_{tss}(D+D_{v})$ is the traffic signal and stop sign status $D_{v}$ in front of the vehicle, with the value of 1 meaning the traffic signal is red or there is a stop sign in front, with the value of 0 meaning the traffic signal is green or there is no stop sign. Variable $D_{sf}$ here indicates the distance to the traffic light or stop sign when no vehicle is in front.

\section{Deterministic Optimal Eco-driving}  \label{se:deter_eco}
The optimal eco-driving control problem is formulated as a nonlinear spatial trajectory optimization problem to minimize vehicle fuel consumption. The cost function $J$ is defined as
\begin{equation}  \label{eq:fuel_cost_J}
\text{minimize}~~J = \int_{0}^{D_{f}} \dot{m}_{fuel} (T_{eng}(D),\omega_{eng}(D) ) \ dD
\end{equation}

The engine torque, wheel braking torque and transmission gear number are chosen as the control variables.
\begin{equation}
u = [T_{eng}(D),T_{brk}(D),N_{gb}(D)]
\end{equation}
The vehicle velocity and traveling/driving time are chosen as the state variables. 
\begin{equation}
x = [v(D), t(D)] 
\end{equation}
\begin{equation}
\frac{dv(D)}{dD} = \frac{a(D)}{v(D)};~\frac{dt(D)}{dD}=\frac{1}{v(D)}
\end{equation}
Subject to the following vehicle physical constraints,
\begin{equation}
\begin{gathered}
T^{min}_{eng} \leq T_{eng}(D) \leq T^{max}_{eng}, \quad \forall \ D \ [0, D_f]  \\
T^{min}_{brk} \leq T_{brk}(D) \leq T^{max}_{brk}, \quad \forall \ D \ [0, D_f]   \\
N_{gb}(D) \in \{ 1,2,3,4,5,6 \}, \quad \forall \ D \ [0, D_f]  \\
v(0) = v(D_{f}) = 0  \\
a^{min} \leq a(D) \leq a^{max}, \quad \forall \ D \ [0, D_f]  \\
v^{min}(D) \leq v(D) \leq v^{max}(D), \quad \forall \ D \ [0, D_f] 
\end{gathered}
\end{equation}
Subject to the following final arrival time and traffic signal passing constraints,
\begin{equation}
t(D_{f}) \leq t_{f}
\end{equation}
\begin{equation}  \label{eq:sig_pass_ct}
c^{i}_{p} \geq c^{i}_{r}
\end{equation}

A key beneficial feature of a spatial trajectory formulation (as opposed to temporal) is that the signal and final destination arrival times do not need to be known \textit{a priori}. A pre-set maximal arrival time constraint $t_{f}$ is imposed on the final state variable $t(D_{f})$, to balance fuel economy and traveling speed. The traffic signal constraint in (\ref{eq:sig_pass_ct}) enforces vehicles to pass through signalized intersections only at green lights.

The above nonlinear optimization problem is solved via dynamic programming adopted from \cite{sundstrom2010implementation}. Detailed formulations are omitted here. Interested readers please refer to \cite{sun2015dynamic,sciarretta2007control}.

\section{Robust Optimal Eco-driving}
\begin{figure}[t]
	\begin{center}
		\includegraphics[trim = 07mm 07mm 07mm 05mm, clip, width=0.50\textwidth]{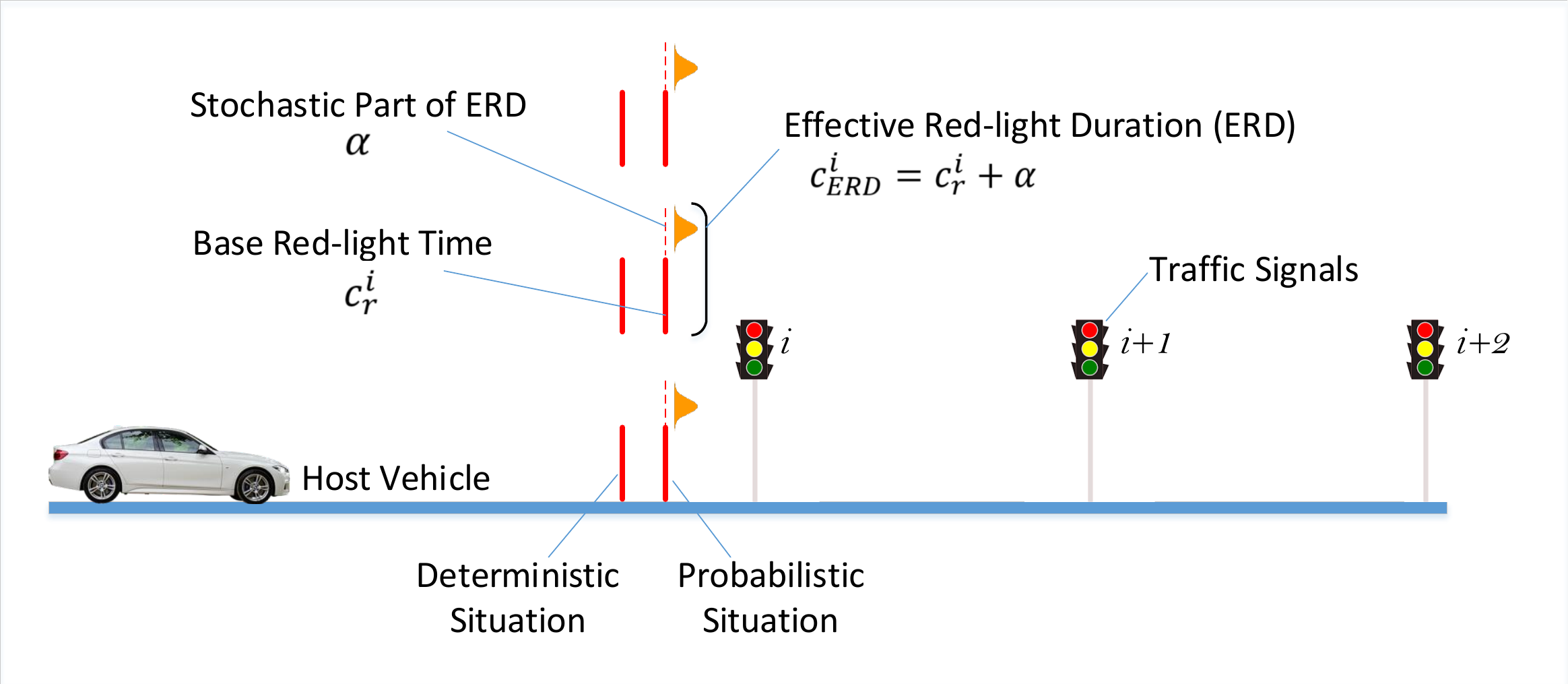}
		\caption{Effective red-light duration, meaning the feasible passing time of a vehicle through an intersection.}
		\label{fig:sto_ecodriving}
	\end{center}
\end{figure}

In Section \ref{se:deter_eco}, it is assumed the SPaT information is deterministic and perfectly known. Mathematically, $c_r^i$ in (\ref{eq:sig_pass_ct}) is known and deterministic. However, as illustrated in Fig. \ref{fig:ERDuncertainty}, the feasible passing time through signalized intersections or crossings is usually uncertain and random. Here, an effective red-light duration (ERD) variable is defined to describe the feasible passing time, denoted as $c^{i}_{ERD}$:
\begin{equation}
c^{i}_{ERD} = c^{i}_{r} + \alpha
\end{equation}
Fig. \ref{fig:sto_ecodriving} exhibits the ERD concept. Parameter $c^{i}_{r}$ is the base red-light duration, which is the minimal red-light time. Random variable $\alpha$ is a stochastic time of delay, caused by signal uncertainties or vehicle waiting queue. In this paper, we assume the total signal cycling-time is not affected by these uncertain factors, meaning $c^{i}_{f}$ is deterministic and known.

Intuitively, $\alpha$ is a random variable over time 0 to $(c^{i}_{f}-c^{i}_{r})$, whose distribution could be (truncated) Poisson, Gaussian, Beta or completely non-parametric. Assume the probability density function of $\alpha$ is $f(\alpha)$. 
Therefore, the traffic signal passing constraint in (\ref{eq:sig_pass_ct}) can be modified to
\begin{equation}  \label{eq:sig_pass_ct_sto}
c^{i}_{p} \geq c^{i}_{ERD} = c^{i}_{r} + \alpha, \qquad \forall \ \alpha 
\end{equation}
However, enforcing the constraint above for all values in the support of $\alpha$ is too restrictive. Consequently, we relax this constraint via chance constraints.

Denote by $\eta$ a required reliability for the subject vehicle to pass through a specific signalized intersection, and $F(\alpha)$ indicates the cumulative distribution function (CDF) of $\alpha$. Equation (\ref{eq:sig_pass_ct_sto}) can be relaxed into the following chance constraint,
\begin{equation}
\text{Pr}(c^{i}_{p} \geq c^{i}_{r} + \alpha ) \geq \eta
\end{equation}
\begin{equation}
\text{Pr}(\alpha \leq c^{i}_{p} - c^{i}_{r} ) = F( c^{i}_{p} - c^{i}_{r} ) \geq \eta
\end{equation}
We assume the CDF $F(\cdot)$ is bijective, and therefore has an inverse function $F^{-1}(\cdot)$. Thus, we can solve for the optimization variable $c^i_p$ to obtain
\begin{equation}
c^{i}_{p} \geq c^{i}_{r} + F^{-1}(\eta)
\end{equation}
Again, $c^{i}_{p}$ is the passing time of subject vehicle through the $i$th intersection in the signal-cycling clock, which is a function of the control and state variables, $c^{i}_{p}(x,u)$.

It should be noted that the real world probability distribution of ERD might vary at different times of day, seasons or locations, and may not be accurately modeled by a parametric distribution. Investigating the actual probability distribution of $\alpha$ via measured data is planned as future work. DP is also used to solve the above robust optimal eco-driving control problem.

\section{Simulation}
The vehicle parameters and engine fuel map used for simulation are extracted from Autonomie \cite{rousseau2014electric}, and summarized in Table \ref{tbl:veh_para}. The maximum and minimal velocity limits are set as 16 and 0m/s, respectively. The acceleration constraints are not activated here, because the engine output torque constraint will restrict the vehicle acceleration within the feasible domain. Three different cases are considered for comparison in the simulation:
\begin{itemize} 
	\item
	Modified IDM. with the human preview-vision distance $D_{v}$ set as 100 meters;
	\item
    Optimal eco-driving with traveling time as cost, denoted as ``Op-time''. Equation (\ref{eq:fuel_cost_J}) is re-formulated as
\begin{equation}
J =  \int_{0}^{D_{f}} t(D) ) \ dD
\end{equation}
	\item
	Optimal eco-driving with fuel consumption as the cost, denoted as ``Op-fuel''.
\end{itemize}

\begin{table}[t]
	\caption{Subject Vehicle Parameters}
	\centering
	\begin{tabular}{p{1.3cm} l p{1.8cm} p{1.3cm}}
		\hline
		Parameter (Unit)    &    Value    &	Parameter (Unit)    &    Value       	\\ [1ex]   
		\hline      
		$m$ (kg)                & 1745	               & $r_{fd}$                                      & 3.51		           \\ [0.6ex]
		$R_{whl}$ (m)        & 0.3413 	           & $\omega^{max}_{eng}$ (rad/s)  	   & 600	             \\ [0.6ex]
		$A$ (m$^2$)          & 2.841                & $T^{max}_{eng}$ (Nm)	                & 240	              \\ [0.6ex]
		$\rho$ (kg/m$^3$)  & 1.1985              & gearbox ratios                & 4.584, 2.964,	\\ [0.5ex]
		$C_{d}$                 &  0.356               &      &  1.912, 1.446,       \\ [0.5ex] 
		$C_{r}$                  & 0.0084, 1.2e-4  &       & 1, 0.74            	\\ [0.5ex]
		\hline                          
	\end{tabular}
	\label{tbl:veh_para} 
\end{table}


\begin{figure*}[!ht]
	\begin{center}
		\includegraphics[trim = 00mm 83mm 00mm 88mm, clip, width=0.65\textwidth]{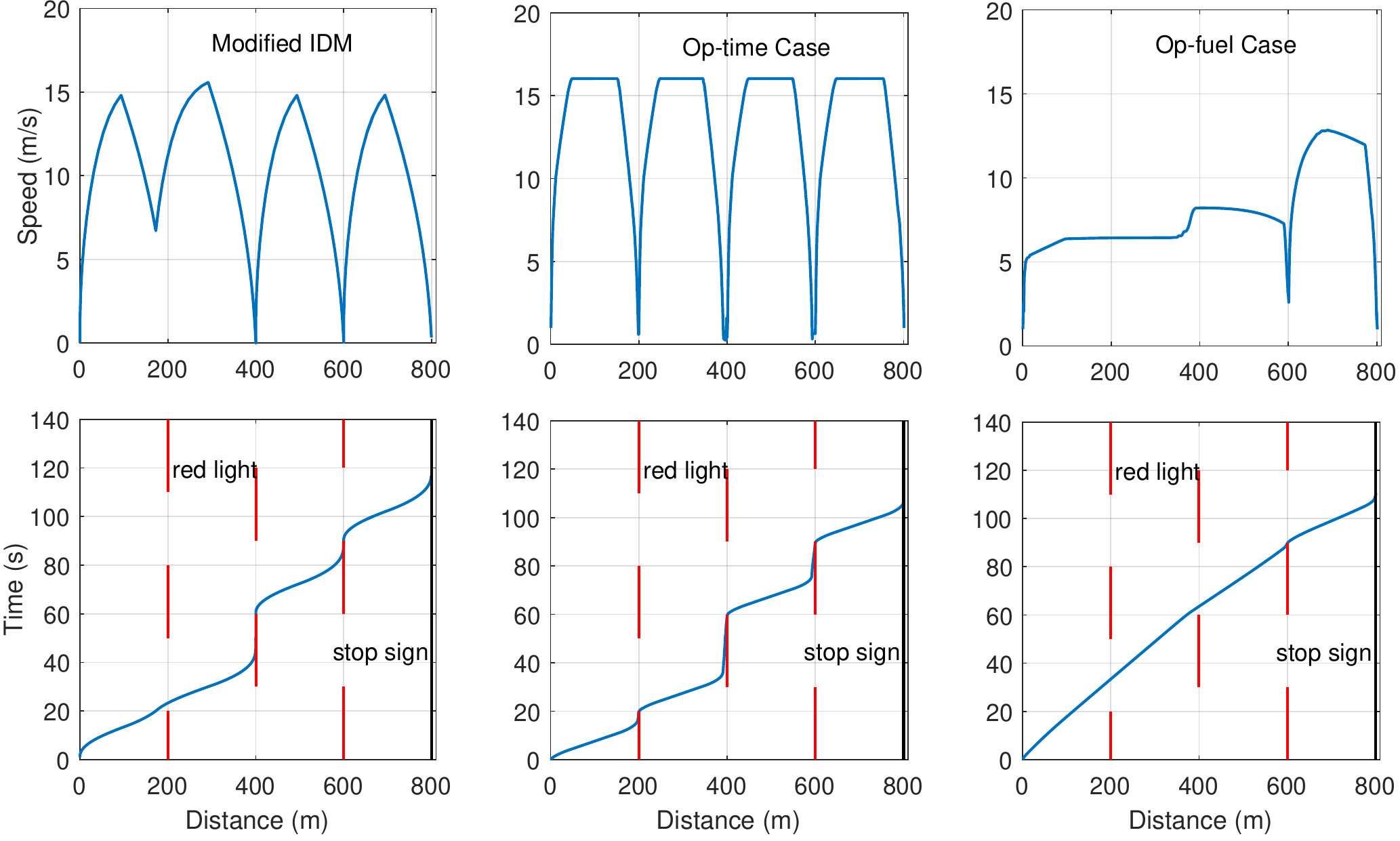}
		\caption{Modified IDM, Op-time and Op-fuel results with deterministic formulation of route 1.}
		\label{fig:route1_res}
	\end{center}
\end{figure*}

\begin{figure*}[!ht]
	\begin{center}
		\includegraphics[trim = 05mm 84mm 03mm 85mm, clip, width=0.657\textwidth]{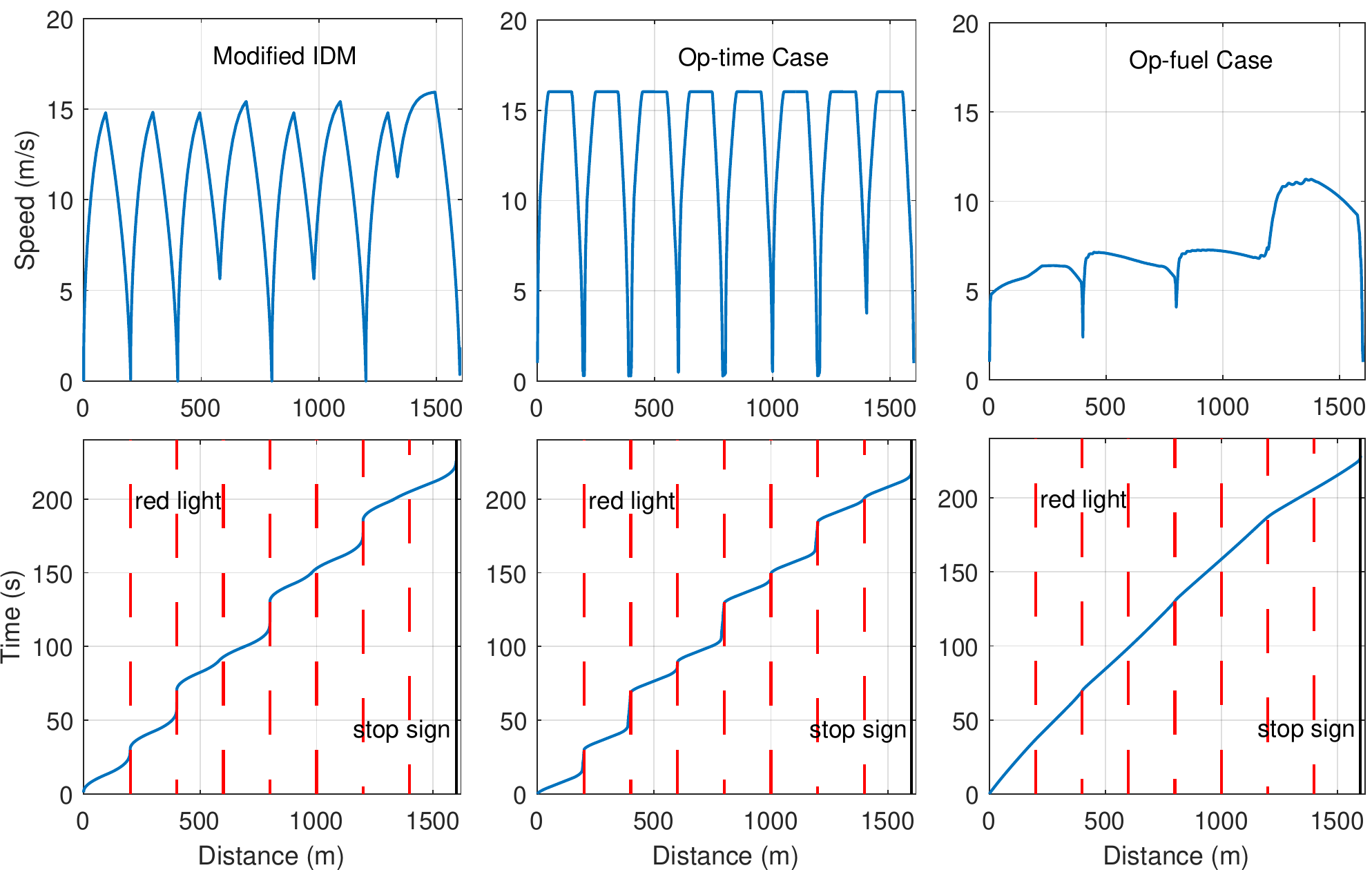}
		\caption{Modified IDM, Op-time and Op-fuel results with deterministic formulation of route 2.}
		\label{fig:route2_res}
	\end{center}
\end{figure*}

\begin{figure}[t]
	\begin{center}
		\includegraphics[trim = 28mm 100mm 31mm 97mm, clip, width=0.46\textwidth]{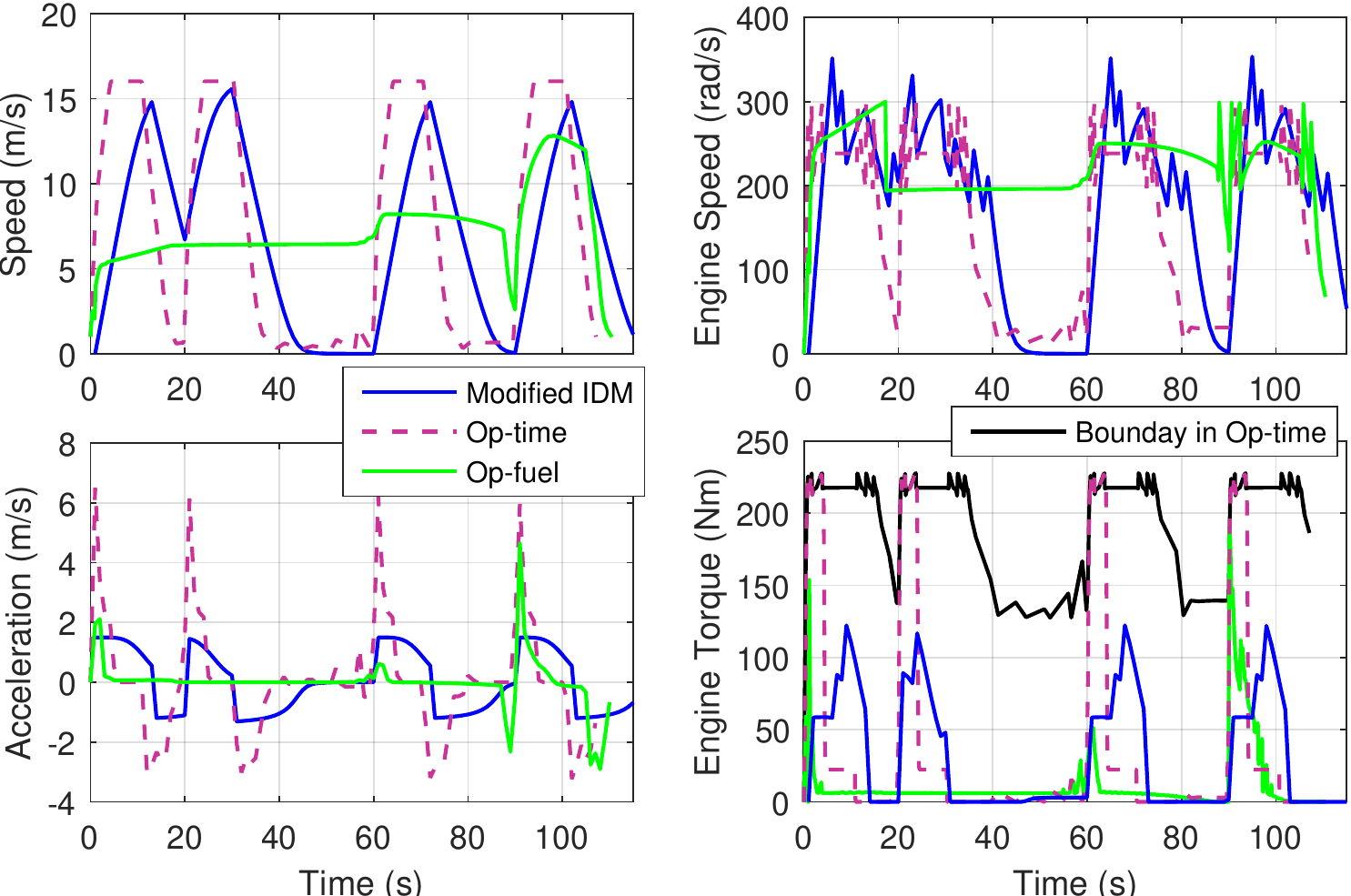}
		\caption{Vehicle velocity, acceleration, engine speed, torque results versus time of route 1.}
		\label{fig:route1_res_detail}
	\end{center}
\end{figure}

\begin{figure}[t]
	\begin{center}
		\includegraphics[trim = 39mm 97mm 45mm 103mm, clip, width=0.39\textwidth]{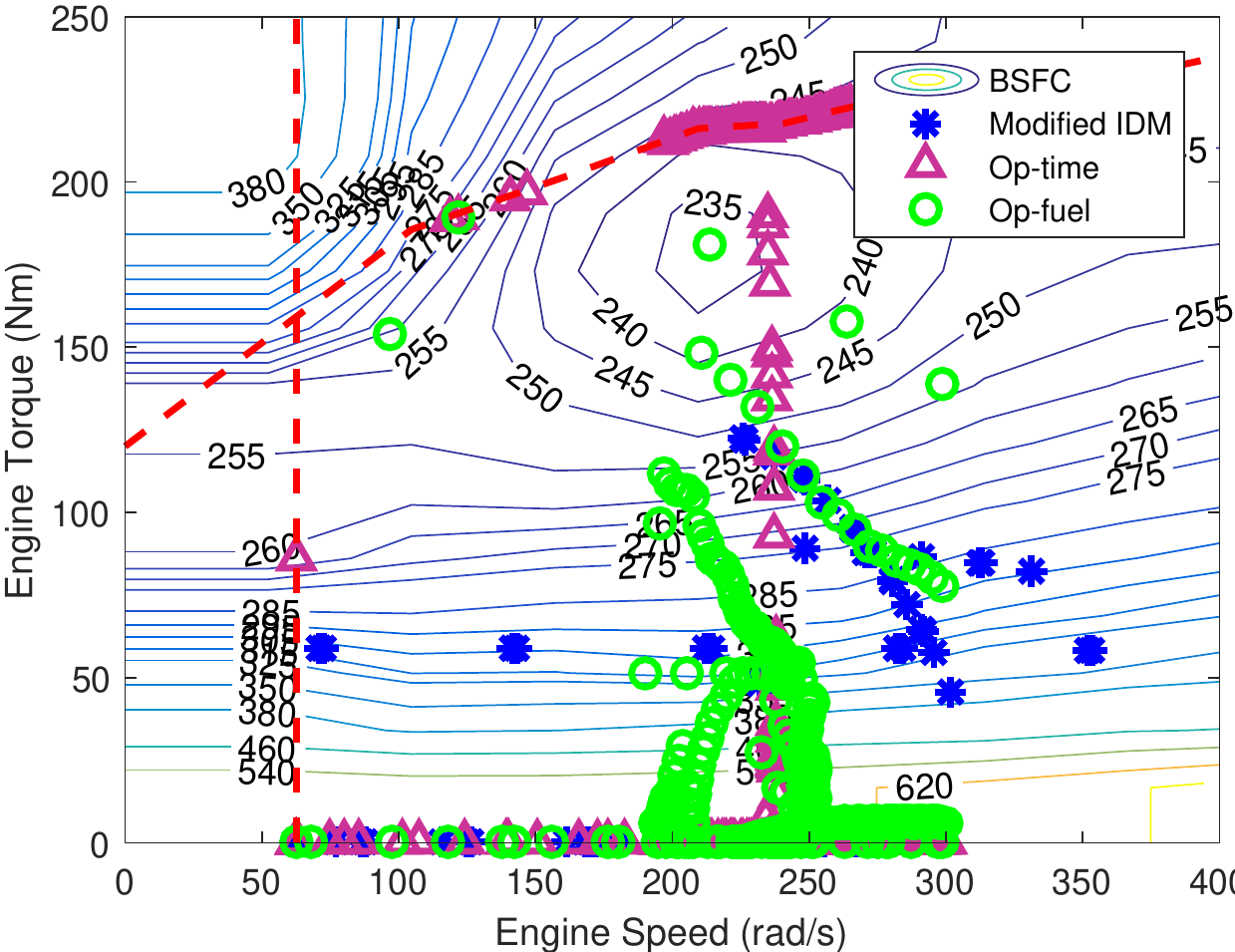}
		\caption{Engine operating points on the BSFC map of route 1.}
		\label{fig:route1_res_eng}
	\end{center}
\end{figure}

\subsection{Deterministic optimal eco-driving}
Two sample driving routes with 3 and 7 signalized intersections, named route 1 and route 2 respectively, are studied in this paper. All of the full cycling periods $c^{i}_{f}$ and red-light durations $c^{i}_{r}$ are intentionally set as 60s and 30s, respectively, for easier analysis of the results. The beginning time $c^{i}_{0}$ is arbitrarily selected between 0 to 30s. However, other realistic selections of the full cycling time and red-light duration can also be incorporated in the proposed optimal eco-driving control strategy. The position and timing information for all the two sample routes are shown in Table \ref{tbl:route_para}.

\begin{table}[h]
	\caption{Position and SPaT information of sample routes}
	\centering
	\begin{tabular}{l | lll | lll | ll}
		\multicolumn{1}{c}{}         &   \multicolumn{3}{c}{Route 1}    &    \multicolumn{3}{c}{Route 2}                 &     \multicolumn{2}{c}{1\&2}           \\ [0.4ex]   
		\hline
		No.    &   Type    &	$D^{i}$ (m)    &    $c^{i}_{0}$*    &   Type    &	$D^{i}$ (m)    &    $c^{i}_{0}$*     &    $c^{i}_{f}$    &    $c^{i}_{r}$      \\ [0.8ex]   
		\hline      
		1       & signal	 & 200        & 10		  &   signal    &  200       &     0            &    60    &   30      \\ [0.5ex]
		2       & signal 	 & 400  	  & 30		  &   signal    &  400       &     20           &    60    &   30      \\ [0.5ex]
		3       & signal     & 600        & 0		  &   signal    &  600       &     0            &    60    &   30      \\ [0.5ex]
		4       & stop       & 800        & --		  &   signal    &  800       &     20           &    60    &   30      \\ [0.5ex]
		5       &            &            &    		  &   signal    & 1000       &     0            &    60    &   30      \\ [0.5ex] 
		6       &            &            &     	  &   signal    &  1200      &     25           &    60    &   30      \\ [0.5ex]
		7       &            &            &     	  &   signal    &  1400      &     10           &    60    &   30      \\ [0.5ex]
		8       &            &            &     	  &   stop      &  1600      &     --           &    60    &   30      \\ [0.5ex]
		\hline                          
		\multicolumn{9}{l}{* Arbitrarily selected values}  \\ [0.2ex]
	\end{tabular}
	\label{tbl:route_para}  
\end{table}

For sample route 1, the vehicle velocity and traveling time results derived from the three driving strategies are plotted in Fig. \ref{fig:route1_res}. It can be seen in the modified IDM approach, the driver started decelerating the vehicle at $D$=100m when the it ‘sees’ a red traffic signal in front. Because of the lack of full SPaT information, the modified IDM is not able to preview the future signal dynamics. About 10 seconds later, it had to switch to accelerate the vehicle again at $D$=170m, as the signal turned to green. This behavior wastes fuel.

At the 2nd traffic signal, the red light blocks the intersection. Modified IDM waits for 20 seconds until the light turns green. A similar scenario happened at the 3rd signalized intersection, but with a shorter waiting time. The vehicle eventually arrived at the stop sign (also its destination) at $t$=117s.

The Op-time strategy accelerates the vehicle whenever possible until the velocity hits the maximal boundary. Under this strategy, the vehicle commonly meets red lights. In the 800-meter-long route assumed in this section, the vehicle cumulatively waited for 50s at all three intersections. The vehicle arrived at the final stop sign at $t$=107s, which is the shortest time among the three driving strategies.

In the Op-fuel case, $t_{f}$ is set as 115s in order to make sure the vehicle arrives to the destination at the same time scale as the other two cases. The Op-fuel strategy smoothly passes through the first two intersections, by adjusting the velocity between 7 and 9m/s. A deeper deceleration happened just before driving through the third signalized intersection ($D$=600m) to wait for a green light. After that, the vehicle velocity restored to about 12m/s to ensure it can arrive the final destination within the time limit. The total driving time in the Op-fuel case is $t$=110.5s, which is 3.5s longer than the Op-time case.

The vehicle velocity and traveling time results for route 2 are shown in Fig. \ref{fig:route2_res}, where similar trends are observed. The modified IDM avoided complete stops at 3 signalized intersections out of 7, with a final arrival time of 226s. This indicates that even without any future information of the traffic signals, the vehicle can still catch green light at normal driving pattern. However, in the Op-time case, the vehicle encountered 6 red lights in order to minimize the arrival time. The final arrival time is 217.9s, which is about 8s (3.5\%) smaller that the modified IDM. As expected, the Op-fuel controller refused to aggressively accelerate the vehicle, and crossed most of the intersections at lower speeds without any complete stops. Eventually, the Op-fuel arrived at the destination at $t$=228.5s.

\begin{figure*}[t]
	\begin{center}
		\includegraphics[trim = 40mm 2mm 51mm 6mm, clip, width=0.93\textwidth]{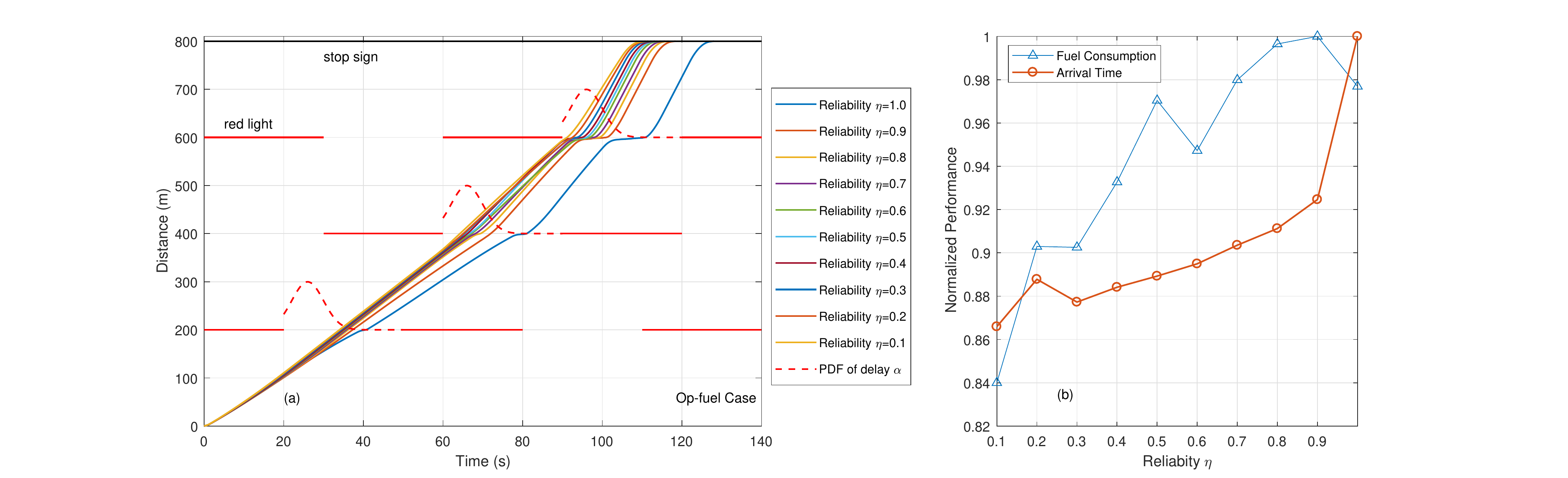}
		\caption{(a) Vehicle traveling trajectories of Op-fuel with different chance reliability $\eta$ used along route 1 under the moderate traffic scenario. (b) Normalized fuel consumption and arrival time results of Op-fuel control with different reliabilities enforced along route 1.}
		\label{fig:sto_route1}
	\end{center}
\end{figure*}

The vehicle velocity, acceleration, engine speed, and engine torque trajectories for route 1 are shown in Fig. \ref{fig:route1_res_detail}. The engine speed is restricted between 200 and 300 rad/s in the Op-fuel case, and the engine torque is relatively smaller than the other two approaches, forming a milder driving style.

The operating points on the brake specific fuel consumption (BSFC) map for route 1 are shown in Fig. \ref{fig:route1_res_eng}. The engine operating point results in Fig. \ref{fig:route1_res_eng} may seem counter-intuitive, but are actually very interesting. The engine operating points from the Op-time case are located more in the high efficiency area of the BSFC map (the lower value the better), compared with the Op-fuel and modified IDM approaches. This is usually preferable in the operation of engines, and often results in better fuel economy. However, the simulated fuel consumption in the Op-time case is, in fact, the highest one and much higher than the Op-fuel result. The main reason is that although the average engine fuel efficiency in the Op-fuel approach is lower, its total engine power requirement is much less. Thus, better fuel economy is achieved with modified IDM and Op-fuel.

The arrival time, average BSFC value, and total engine fuel consumption results of route 1 and 2 are reported in Table \ref{tbl:route_res}. The average BSFC of Op-fuel is 16.9-21.2\% higher than that of Op-time, while the overall fuel consumption is 51.9-59.5\% less. This significant fuel economy improvement is achieved by sacrificing 2.8-4.9\% of the arrival time, which is trivial in daily driving.

\begin{table}[h]
	\caption{Arrival time, average BSFC and fuel consumption results of route 1 and 2 with deterministic SPaT}
	\centering
	\begin{tabular}{l | llll}
		\hline
		Route  &  Method          &   $t(D_{f})$ (s)    &   $\text{B}_{avg}^{\diamond}$ (g/kWh)     &   Fuel (g)         \\ [0.8ex]   
		\hline      
                     & Modified IDM	    & 117                 & 478.26		      &   88.24           \\ [0.5ex]
		800m         & Op-time 	        & 107  	              & 485.94            &   91.49           \\ [0.5ex]
		3 lights     & Op-fuel          & 110                 & 568.19            &   43.95           \\ [0.5ex]
		             & Change*          & +2.8\%              & +16.9\%           &   -51.9\%         \\ [0.5ex]
		\hline
		             & Modified IDM      & 226                & 477.31            &   172.84          \\ [0.5ex] 
		1600m        & Op-time           &  217.9             & 486.62	          &   182.15          \\ [0.5ex]
		7 lights     & Op-fuel           & 228.52             & 586.53            &   73.79           \\ [0.5ex]
		             & Change*           & +4.9\%             & +21.2\%           &   -59.5\%          \\ [0.5ex]
		\hline                          
		\multicolumn{5}{l}{* Indicates the performance change of Op-fuel compared with Op-time;}           \\ [0.2ex]
		\multicolumn{5}{l}{$\diamond$ $\text{B}_{avg}$ is the average BSFC value.}           \\ [0.2ex]
	\end{tabular}
	\label{tbl:route_res}  
\end{table}

\begin{figure}[t]
	\begin{center}
		\includegraphics[trim = 35mm 85mm 34mm 83mm, clip, width=0.42\textwidth]{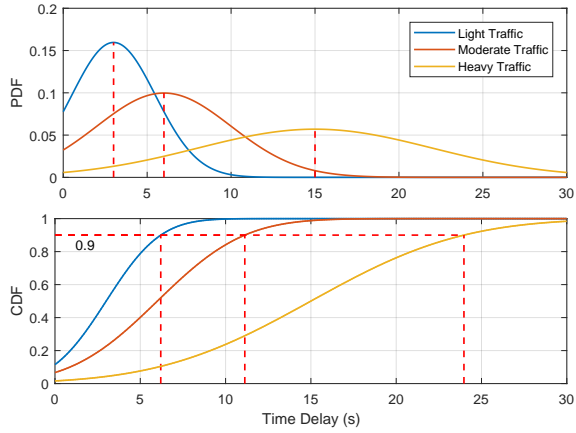}
		\caption{Probability density function and cumulative density function of $\alpha$.}
		\label{fig:pro_curve}
	\end{center}
\end{figure}


\begin{figure}[t]
\begin{center}
	\includegraphics[trim = 37mm 78mm 43mm 78mm, clip, width=0.4\textwidth]{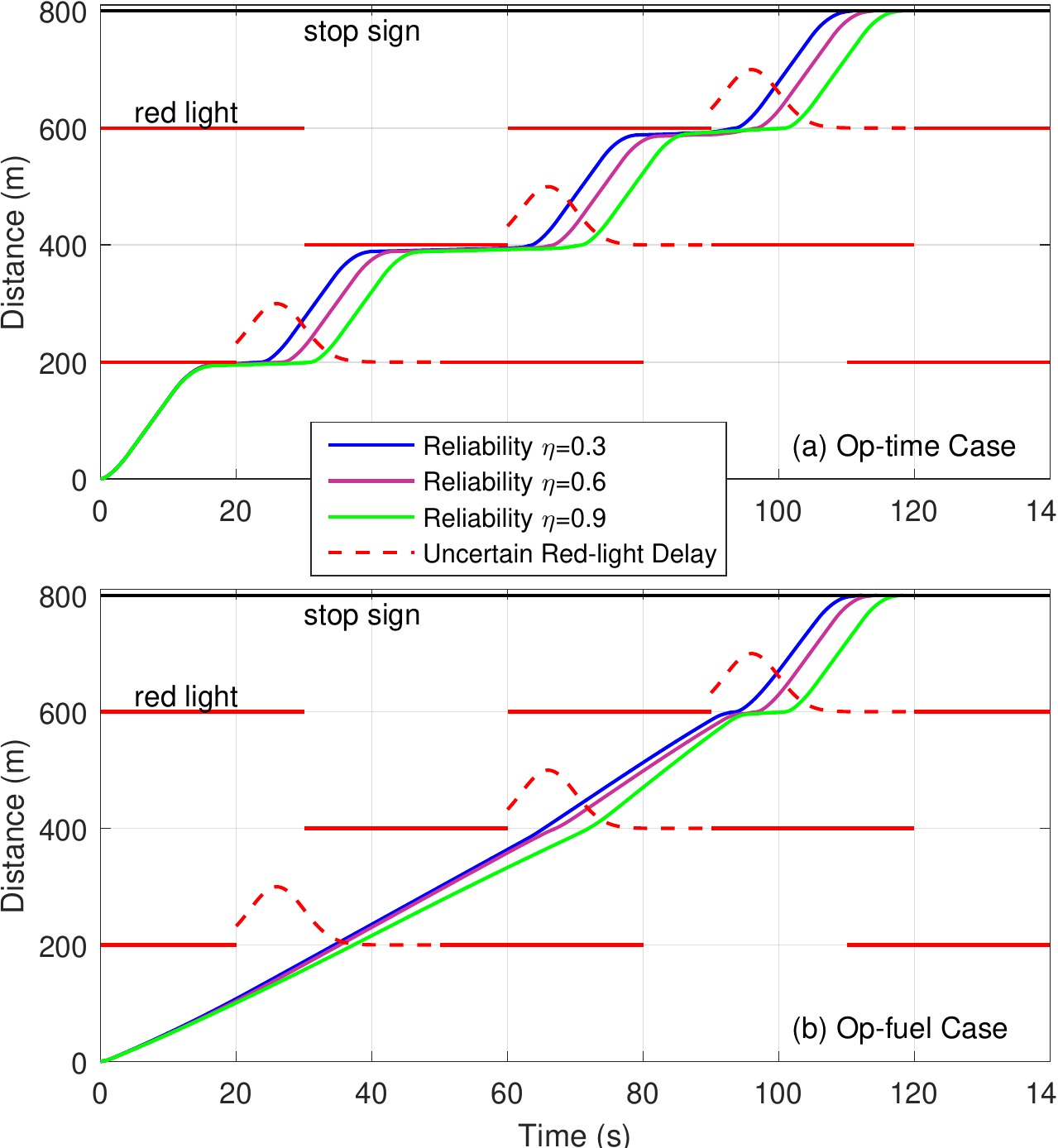}
	\caption{Comparison of robust Op-time and robust Op-fuel control results when $\eta$=0.3,0.6,0.9, respectively.}
	\label{fig:sto_resall}
\end{center}
\end{figure}

\subsection{Robust optimal eco-driving}
In this section, we assume $\alpha$ is a truncated Gaussian random variable for illustrative purposes, with the PDF and CDF drawn in Fig. \ref{fig:pro_curve}. The true probability distribution for the ERD is generally unknown. Future work will focus on this challenge.

Three possible scenarios are shown in Fig. \ref{fig:pro_curve}: light, moderate and heavy traffic situations. When traffic is light, the high-probability values of $\alpha$ are around 0 to 6s, indicating small time delays. For heavy traffic scenarios (for example, at rush hours), the average $\alpha$ value is much larger at around 15s. It is even possible that the vehicle would have to wait for the next green signal, which is normal in real life. The moderate traffic situation is adopted for simulation, where
\begin{equation}
\alpha \sim \mathcal{N}_{[0,30]}(6,16) \in [0,30]
\end{equation}
The vehicle traveling trajectories of Op-fuel control along route 1, with various chance reliability $\eta$ are illustrated in Fig. \ref{fig:sto_route1}.

As can be seen in Fig. \ref{fig:sto_route1}(a), the arrival time increases as reliability parameter $\eta$ increases. For $\eta = 0.1$, the optimal velocity trajectory has very little robustness to delay beyond the base red-light time. This can be clearly seen from the vehicle trajectories at the 600-meter traffic signal, where the vehicle passed the intersection immediately after the light switched to green. As $\eta$ increases, the passing time gradually increases as the solution becomes more cautious to delayed ERD. Yet at the 200-meter and 400-meter signalized intersections, the originally planned passing time has already avoided most of the possible delays. Thus, the trajectories with chance reliability $\eta$ from 0.1 to 0.8 are quite identical. 

Fig. Fig. \ref{fig:sto_route1}(b) shows the normalized fuel consumption and arrival time results of Op-fuel with different reliabilities enforced along route 1. When $\eta$ is 0.9, the fuel consumption increased 16\% compared with $\eta = 0.1$. The arrival time change is much smaller, with an increase of only 5\%. At the $\eta$=1.0 point, the arrival time is raised by nearly 8\%, yet the fuel consumption decreased by 2\%. Fig. \ref{fig:sto_resall} is the comparison of Op-time and Op-fuel results under the robust optimal eco-driving control formulation, when $\eta$=0.3,0.6,0.9 respectively. The main difference is that the Op-time control is less tolerant to the delay uncertainty, and thus its performance decline is more severe than the Op-fuel control.

Table \ref{tbl:route_res_all} summaries the arrival time, average BSFC and fuel consumption results of robust optimal eco-driving and modified IDM with both driving route 1 and 2. Comparing deterministic control with robust control across the three methods, we find that the final arrival time grows as the reliability $\eta$ increases. However, the fuel consumption increase is not as significant. Obviously, the fuel consumption of Op-fuel is less than that of Op-time and modified IDM. Its arrival time is slightly longer than Op-time, but mostly shorter than modified IDM.

\begin{table}[h]
	\caption{Arrival time, average BSFC and fuel consumption results of route 1 and 2 with stochastic SPaT}
	\centering
	\begin{tabular}{l | ll | lll}
		\hline
		Route  &  Method   & Type       &   $t(D_{f})$    &   $\text{B}_{avg}^{\diamond}$     &   Fuel        \\ [0.8ex]   
		\hline      
		               &  Modified  & Deterministic              & 117            & 478.26		  & 88.24         \\ [0.5ex]
		               &  IDM	    & Robust $\eta$=0.3 	     & 121  	      & 479.33        & 89.48         \\ [0.5ex]
		               &            & Robust $\eta$=0.6          & 125            & 478.10        & 91.09          \\ [0.5ex]
		               &            & Robust $\eta$=0.9          & 129            & 481.93        & 90.63          \\ [0.5ex]
	                   &  Op-time  	& Deterministic              & 107.2          & 485.94        & 91.49         \\ [0.5ex] 
		800m           &            & Robust $\eta$=0.3          & 111.1          & 485.52	      & 92.11         \\ [0.5ex]
		3 lights       &    	    & Robust $\eta$=0.6          & 114.0          & 485.31        & 92.22          \\ [0.5ex]
		               &   	        & Robust $\eta$=0.9          & 118.2          & 485.38        & 91.78          \\ [0.5ex]
				       &  Op-fuel   & Deterministic              & 110.4          & 568.19        & 43.95         \\ [0.5ex] 
		               &    	    & Robust $\eta$=0.3          & 112.3          & 557.59	      & 48.14         \\ [0.5ex]
		               &    	    & Robust $\eta$=0.6          & 114.6          & 562.23        & 50.43          \\ [0.5ex]
		               &   	        & Robust $\eta$=0.9          & 118.4          & 567.61        & 53.16          \\ [0.5ex]
		\hline         
				       &  Modified 	  & Deterministic            & 226            & 477.31        & 172.84         \\ [0.5ex] 
		               &  IDM  	      & Robust $\eta$=0.3        & 229            & 478.90	      & 172.72         \\ [0.5ex]
		               &   	          & Robust $\eta$=0.6        & 233            & 480.97        & 172.71          \\ [0.5ex]
		               &   	          & Robust $\eta$=0.9        & 237            & 482.98        & 172.70          \\ [0.5ex]
				       &  Op-time 	  & Deterministic            & 217.9          & 486.62        & 182.15         \\ [0.5ex] 
		1600m          &    	      & Robust $\eta$=0.3        & 221.5          & 486.66	      & 181.24         \\ [0.5ex]
		7 lights       &    	      & Robust $\eta$=0.6        & 224.6          & 487.08        & 181.80          \\ [0.5ex]
		               &   	          & Robust $\eta$=0.9        & 227.7          & 485.87        & 183.18          \\ [0.5ex]
				       &  Op-fuel 	  & Deterministic            & 228.52         & 586.53        & 73.79        \\ [0.5ex] 
		               &    	      & Robust $\eta$=0.3        & 229.6          & 579.84	      & 72.66        \\ [0.5ex]
		               &    	      & Robust $\eta$=0.6        & 230.8          & 576.53        & 72.98         \\ [0.5ex]
		               &   	          & Robust $\eta$=0.9        & 232.4          & 564.48        & 72.85         \\ [0.5ex]
		\hline         
		\multicolumn{5}{l}{$\diamond$ $\text{B}_{avg}$ is the average BSFC value.}           \\ [0.2ex]
	\end{tabular}
	\label{tbl:route_res_all}  
\end{table}

\section{Conclusions and Future Work}

This paper proposes a novel robust optimal eco-driving control strategy to solve the vehicle velocity planning problem with multiple signalized intersections, based on a spatial optimization formulation. The requirement for prior knowledge of the destination arrival time is eliminated. We propose a novel traffic signal modeling approach. Effective red-light duration (ERD) is proposed to capture the random feasible passing time at signalized intersections. The optimal control problem is solved via dynamic programming (DP). Simulation results indicate that the developed optimal eco-driving strategy is able to reduce fuel consumption by approximately 50-57\%, while maintaining the arrival time at the same level compared with the modified intelligent driver model. The controller robustness to signal timing uncertainty is greatly improved with slight sacrifices to vehicle fuel economy.

Future work includes real-world traffic SPaT probability distribution study. We also plan to develop methods to reduce the optimal eco-driving control computation complexity.




\ifCLASSOPTIONcaptionsoff
  \newpage
\fi

\bibliographystyle{IEEEtran}
\bibliography{IEEEabrv,Robust_op_eco}

\end{document}